\documentclass{arXiv}
\usepackage{amssymb, amsmath}

\newtheorem{theorem}{Theorem}[section] 
\newtheorem{lemma}[theorem]{Lemma}     
\newtheorem{corollary}[theorem]{Corollary}

\newnumbered{assertion}{Assertion}    
\newnumbered{conjecture}{Conjecture}  
\newnumbered{definition}[theorem]{Definition}
\newnumbered{hypothesis}{Hypothesis}
\newnumbered{remark}[theorem]{Remark}
\newnumbered{note}{Note}
\newnumbered{observation}{Observation}
\newnumbered{problem}{Problem}
\newnumbered{question}{Question}
\newnumbered{algorithm}{Algorithm}
\newnumbered{example}{Example}
\newunnumbered{notation}{Notation} 

\providecommand{\R}{\mathbb{R}}
\providecommand{\N}{\mathbb{N}}
\providecommand{\s}{\mathbb{S}}

\providecommand{\h}{\mathcal{H}}
\providecommand{\lr}{\langle\cdot,\cdot\rangle}

\providecommand{\pp}{\mathcal{P}_2({\R^d})}
\providecommand{\ppz}{\mathcal{P}_{2,h}({X})}
\providecommand{\ppo}{\mathcal{P}_{2,0}({X})}

\providecommand{\px}{\mathcal{P}_2(X)}
\providecommand{\ppx}{\mathcal{P}_{\!p}(X)}
\providecommand{\py}{\mathcal{P}_2(Y)}

\providecommand{\mt}{\widetilde{\mu}}
\providecommand{\mmp}{\mu^{\prime}}
\providecommand{\nt}{\widetilde{\nu}}

\providecommand{\x}{x^{\prime}}
\providecommand{\y}{y^{\prime}}
\providecommand{\zp}{{z^{\prime}}}
\providecommand{\z}{{h^{\prime}}}
\providecommand{\tp}{{s^{\prime}}}

\providecommand{\vo}{\varphi_\z}
\providecommand{\vz}{\varphi_h}
\providecommand{\vzs}{(\vz)_ {\sharp}}

\providecommand{\tsl}[1]{\textsl{#1}}

\title[Cone structure of $L^2$-Wasserstein spaces]
 {Cone structure of $L^2$-Wasserstein spaces} 
\author[A. Takatsu \and T. Yokota]{Asuka Takatsu \and Takumi Yokota}
\classno{60D05 (primary), 51K10, 58B20}
\extraline{Keywords:Wasserstein space, cone structure, rank, splitting theorem}
 
\begin{document}
\maketitle
\begin{abstract}
The purpose of this paper is to understand the geometric structure of 
the $L^2$-Wasserstein space $\pp$ over the Euclidean space.
For this sake, we focus on its cone structure.
One of our main results is that 
the $L^2$-Wasserstein space over a Polish space has a cone structure if and only if so does the underlying space.
In particular, $\pp$ turns out to have a cone structure.
It is also shown that $\pp$ splits $\R^d$ isometrically but not $\R^{d+1}$.
\end{abstract}

\section{Introduction}
Let $(X, d)$ be a Polish space, i.e., a complete separable metric space.
The $L^p$-Wasserstein space over $(X, d)$ is
the set $\ppx$ of all Borel probability measures on $X$ with finite $p$-th moment,
endowed with the so-called $L^p$-Wasserstein distance denoted by $W_p$. 
The definition of $W_p$,
which is recalled in the next section, has its root in the optimal transport theory.

Since the theory was born,
contributions have been made by a number of authors.
A milestone was the discovery made by Otto~\cite{ot}
that the solutions of porous medium equations can be regarded as gradient flows
on the $L^2$-Wasserstein space $(\pp, W_2)$ (cf.~\cite{jko}).
This was done by introducing a formal Riemannian structure to $\pp$
whose induced distance coincides with the $L^2$-Wasserstein distance.
In addition, he performed formal calculation to demonstrate that
the sectional curvature of $\pp$ with respect to his formal structure is everywhere non-negative.

Later, the non-negativity of the curvature was justified by showing that
\tsl{ $(\px, W_2)$ is an Alexandrov space of non-negative curvature if and
only if so is the underlying space $(X, d)$} 
(e.g. Sturm~\cite{st}).
This fact suggests a close relationship between the geometry of $\px$
and that of the underlying space $X$.

The purpose of the present paper is to contribute to a better understanding
of the geometric structure of the $L^2$-Wasserstein space
over the Euclidean space $\R^d$.
Among various special structures of $\R^d$,
we shall focus on its cone structure.
The Euclidean space is naturally isometric to the (Euclidean) cone
of its unit sphere with the angle metric $\angle$
(see Definition~\ref{def;cone} below).
Now the main theorem of this paper is formulated as follows:    
\begin{theorem}\label{ippan}
Let $(X, d)$ be a Polish space.
Then its $L^2$-Wasserstein space $(\px, W_2)$ 
has a cone structure 
if and only if
so does the underlying space $(X, d)$.
Furthermore, if this is the case,
$\px$ is non-branching at the vertex
if and only if
so is the underlying space $X$.
\end{theorem}

In the statement above,
we say that a metric space $(X, d)$ is 
\tsl{non-branching} at a point $x$ in $X$
if for any points $u, v, w$ in $X$,
\begin{equation}\label{non-branch}
\begin{split}
d(u, v) &= d(u, x)+d(x, v)\\
&= d(u, x)+d(x, w) = d(u, w)
\end{split}
\quad \implies v=w.
\end{equation}

When $(X, d)$ and $x$ happen to be 
a cone over a metric space $(\Sigma, \angle)$ 
and its vertex respectively,
this definition is equivalent to 
that $(\Sigma, \angle)$ satisfies for any $\xi$ in $\Sigma$,
$$\text{the antipodal set } \{\xi^{\prime} \in \Sigma \ |\ \angle(\xi^{\prime}, \xi) \ge \pi\} 
\text{ consists of at most one point}.$$
Remember that any Hilbert space is a cone 
which is non-branching at every point.

The main ingredient of the proof of Theorem~\ref{ippan} is
the analysis of the behavior of $L^2$-Wasserstein geodesics,
particularly those passing through Dirac measures.
As a consequence,
it will be shown that when an $L^2$-Wasserstein space has a cone structure,
its vertex must be a Dirac measure (Corollary~\ref{vertex}).
This observation plays a crucial roll in the proof of the ``only if'' part,
as well as in the proofs of the corollaries stated below.

In the previous preprint~\cite{asuka} of the first author,
she found out that the space of all Gaussian measures,
with the metric induced from $(\pp, W_2)$, has a cone structure.
The structure of the base space are also discussed in~\cite{asuka}.
By restricting to such a subset of $\pp$,
which is a finite dimensional Riemannian manifold,
she also gave a justification to Otto's calculation of the sectional curvature.

We continue our study on the geometric structure of $L^2$-Wasserstein spaces.
It is clear that $\pp$ contains an isometric copy of the underlying space $\R^d$.
We shall obtain more information.

At first, following Kloeckner~\cite{Kl},
we define the \tsl{rank} of a metric space
as the supremum of the dimensions of which
Euclidean spaces can be isometrically embedded into it 
(cf.~Foertsch--Schroeder~\cite{FS}).
We know that the rank of the $L^p$-Wasserstein space
is not less than that of the underlying space.
One of the main results of a recent preprint~\cite{Kl} is
the coincidence of the ranks of $\R^d$ and $\pp$.

\begin{corollary}\label{rank}
Let $(X, d)$ be a Polish space which has a cone structure and 
is non-branching at the vertex.
Consider the tower $\{X_i\}_{i=0}^\infty$ of Polish spaces
constructed by letting $X_0 =X$ and 
$X_{i+1}$ be the $L^2$-Wasserstein space over $X_i$.
Then all of the ranks of $X_i$'s are equal to 
that of the underlying space $X=X_0$.
\end{corollary}

Next we establish a splitting theorem for $L^2$-Wasserstein spaces.
A metric space $(X, d)$ is said to split a Hilbert space $(\h,\lr)$
when $(X, d)$ is isometric to the direct product of 
$(\h, \lr)$ and some metric space.
We prove the following theorem for general Polish spaces.
\begin{theorem}\label{split}
If a Polish space $(X, d)$ splits a separable Hilbert space $(\h, \lr)$,
then its $L^2$-Wasserstein space $(\px, W_2)$ splits $(\h, \lr)$ as well.
\end{theorem}

As a by-product of our work,
we obtain a partial converse to the previous theorem.
\begin{corollary}[(of Theorem~\ref{ippan})]\label{tilps}
Let $(X, d)$ be a Polish space.
Suppose that its $L^2$-Wasserstein space $(\px, W_2)$
splits a separable Hilbert space $(\h, \lr)$.
If either
\begin{enumerate}
\item $(X, d)$ is non-branching at every point, or
\item $(X, d)$ has a cone structure (and hence so does $(\px, W_2)$),
\end{enumerate}
then $(X, d)$ also splits $(\h, \lr)$.
\end{corollary}

Finally,
we summarize our results applied to the Euclidean space as follows.
\begin{corollary}
The $L^2$-Wasserstein space $\pp$ over the Euclidean space 
$\R^d$ has a cone structure,
and is isometric to the direct product of 
$\R^d$ and a certain Polish space whose rank is zero.
\end{corollary}
As far as the authors know,
this is a result which is not available in the literature.
(While preparing this paper,
the authors came across the paper by Carlen--Gangbo~\cite{cg}
which is closely related to our this work.
See Remark~\ref{cg} below.)

The organization of this paper is as follows:
In the next section, we recall the necessary definitions
from metric and Wasserstein geometry.
In Section 3, we describe the proof of Theorem~\ref{ippan}.
Section 4 is devoted to the proofs of
Theorem~\ref{split} and Corollaries~\ref{rank} and~\ref{tilps}.

\section*{Acknowledgements}
The first author is grateful to Sumio Yamada for many valuable suggestions and encouragement.
She would also like to thank Takashi Shioya for his suggestions and comments,
and Kazumasa Kuwada for his remarks, 
one of which provided a starting point for this work.
Both authors are indebted to Shin-ichi Ohta and Masayoshi Watanabe
for their comments.

The authors also acknowledge the support of 
the Research Fellowships of 
the Japan Society for the Promotion of
Science for Young Scientists.

\section{Preliminaries}
\subsection{Background on metric spaces}
In this subsection, 
we summarize some definitions on the geometry of metric spaces.
For further details, 
we refer to~\cite{BBI} and~\cite{BGP}.
Let $(X, d)$ be a metric space.
\begin{definition}
Let $x,y,z$ be three distinct points in $X$.
We denote by $\tilde{\angle} x y z$ 
the \tsl{comparison angle} of $\angle x y z$, which is defined by 
\[
  \tilde{\angle} x y z
=\arccos
 \frac{d(x,y)^2+d(y,z)^2-d(z,x)^2}
  {2d(x,y)d(y,z)}.
\]
\end{definition}
\begin{definition}
Let $\gamma:[0,\varepsilon)\to X$ 
and 
$\sigma:[0,\varepsilon)\to X$ 
be two paths in $X$
starting at the same point $x$.
We define the \tsl{angle} $\angle_x (\gamma,\sigma)$ 
between $\gamma$ and $\sigma$ as
\[
  \angle_x (\gamma,\sigma)
=\lim_{s,t \searrow 0} 
 \tilde{\angle} \gamma(s)x\sigma(t),
\]
if the limit exists.
\end{definition}
We briefly discuss the tangent cone of $X$.
Fix a point $x$ in $X$. 
We assume that $\angle_x (\gamma,\sigma)$ always exists 
for any two geodesics.
A \tsl{geodesic} is a constant speed curve 
whose length is equal to the distance 
between its endpoints.
We define $\Sigma^{\prime}_x$ as the set of all geodesics 
starting at $x$ 
equipped with an equivalence relation ${}\parallel{}$, 
where   
$\gamma \parallel \sigma$ 
means $\angle_x (\gamma,\sigma)=0$.
The angle $\angle _x$ is independent of the choices of 
$\gamma$ and $\sigma$ in their equivalence classes.
Then $\angle _x$ is a natural distance function 
on $\Sigma^{\prime}_x$. 
We define the \tsl{space of directions} 
$(\Sigma_x, \angle_x)$ at $x$ as the metric 
completion of 
$(\Sigma^{\prime}_x, \angle_x)$.
The \tsl{tangent cone} $(K_x, d_x)$ at $x$
is, by definition, the cone over $(\Sigma_x, \angle_x)$.
\begin{definition}\label{def;cone}
The \tsl{cone} over a metric space $(\Sigma, \angle)$ is the quotient space 
$C(\Sigma)=\Sigma \times [0,\infty)/\sim{}$, 
where the equivalence relation $\sim$ is defined by
$(\xi,s) \sim (\eta,t)$ 
if and only if $s=t=0$.
We call the equivalence class of $(\cdot,0)$ and $\Sigma$ 
the \tsl{vertex} and the \tsl{base space}, respectively.
The distance $d_C$ on the cone is defined by
\[
 d_C((\xi,s), (\eta,t))
=\sqrt{s^2+t^2-2st \cos (\min\{\angle(\xi,\eta),\pi\})}. 
\]
A metric space is said to have a cone structure 
when it is isometric to some cone. 
\end{definition}
It is trivial that 
$(\Sigma, \angle)$ is a Polish space 
if and only if so is $(C(\Sigma), d_C)$. 
%
\subsection{$L^p$-Wasserstein spaces}
In this subsection, 
we review $L^p$-Wasserstein spaces (see~\cite{V},~\cite{Vi}).
Let $(X, d)$ be a Polish space.
Given two Borel probability measures $\mu$ and $\nu$ on $X$,
a \tsl{transport plan} $\pi$ between $\mu$ and $\nu$ 
is a Borel probability measure on $X \times X$ 
with marginals $\mu$ and $\nu$, that is, 
\[
\pi [ A \times X ]=\mu [A],
\quad \pi [ X \times A]= \nu [A]
\]
for all Borel sets A in $X$.
We denote by $\Pi(\mu,\nu)$ 
the set of transport plans between $\mu$ and $\nu$.
\begin{definition}
For any two Borel probability measures $\mu$ and $\nu$ on $X$, 
the \tsl{$L^p$-Wasserstein distance} between $\mu$ and $\nu$
is defined by
\[
W_p(\mu ,\nu ) = \left(
                 \inf _{\pi \in \Pi (\mu, \nu)}
                 \int _{X \times X}d(x,y)^p d \pi (x,y) \right) ^{\frac1p}.
\]
\end{definition}
In general, this does not define a distance 
function on the set of 
all Borel probability measures 
because $W_p(\mu,\nu)$ might take the value $\infty$ 
when one of the measures has infinite $p$-th moments.
Henceforth we restrict $W_p$ to 
the set $\ppx$ of all Borel probability measures 
whose $p$-th moments are finite.
Then $W_p$ defines a distance on $\ppx$ 
for $p$ in $[1,\infty)$ and 
we call the pair $(\ppx,W_p)$ 
the \tsl{$L^p$-Wasserstein space} over $(X, d)$.
A transport plan in $\Pi(\mu,\nu)$ 
is said to be \tsl{optimal} 
if it achieves the distance $W_p(\mu,\nu)$.
An optimal transport plan always exists.
Details can be found in~\cite[Chapter~4]{Vi}.

The underlying space $X$ is isometrically embedded into $\ppx$
by identifying a point $x$ in $X$ 
with the Dirac measure $\delta_x$ in $\ppx$.
In particular, if $\gamma(t)$ is a geodesic in $X$, 
then $\delta_{\gamma(t)}$ is a geodesic in $\ppx$.
This face partly demonstrates that
$L^p$-Wasserstein spaces is often adapted to statements that 
combine weak convergence and geometry of their underlying spaces.
In particular,
we stress that 
$L^p$-Wasserstein space over a Polish space 
is itself a Polish space
(see~\cite[Chapter~6]{Vi}).

Since the $L^2$-Wasserstein space 
has a closer relationship with
the ``Riemannian'' geometry of the underlying space 
than the $L^p$-Wasserstein space 
as mentioned in the introduction, 
we treat especially the case $p=2$ in the rest of the paper.

\section{The proof of the main theorem}
We first prove the ``if'' part of Theorem~\ref{ippan}, 
namely, we show that 
if a Polish space $(X, d)$ is a cone 
over $Y$ then
its $L^2$-Wasserstein space $(\px,W_2)$ also
has a cone structure.
Let $o$ and $\delta=\delta_o$ be 
the vertex of $X$ and
the Dirac measure centered at $o$, respectively.
We need to prove some lemmas 
under the assumptions of Theorem~\ref{ippan}.
\begin{lemma}\label{entyou}
For any non-negative number $s$, 
we define a map 
$\psi_s$ on $X=C(Y)$ by $\psi_s(y,t)=(y,st)$ and
the associated map $\Psi_s$ on $\px$ 
by $\Psi_s(\mu)=(\psi_s)_{\sharp} \mu$. 
Then $\{\Psi_s(\mu)\}_{s \in[0,1]}$ 
is a geodesic from $\delta$ to $\mu$ in $\px$.
\end{lemma}
\begin{proof}
Due to the fact that $\Psi_0(\mu)=\delta$ 
and $\Psi_1(\mu)=\mu$,
we only need to show that
\[ 
W_2(\Psi_s(\mu),\Psi_t(\mu))
\leq |s-t| W_2(\delta,\mu)
\] 
for all $s,t$ in $[0,1]$.
Since $(\psi_s \times \psi_t)_{\sharp}\mu$ 
is a transport plan in $\Pi(\Psi_s(\mu),\Psi_t(\mu))$, 
we have
\begin{align*}
   W_2(\Psi_s(\mu),\Psi_t(\mu))^2 
 &\leq \int_{X \times X} 
  d(x_1,x_2)^2 d (\psi_s \times \psi_t)_{\sharp}\mu(x_1,x_2)\\
 &= \int_{X} 
  d(\psi_s(x),\psi_t(x))^2 d \mu(x) \\
 &=(s-t)^2W_2(\delta,\mu)^2.
\end{align*}
\end{proof}
\begin{lemma}\label{kyori}
For any $\mu,\nu$ in $\px$ and non-negative numbers $s,t$, 
we have
\[
  W_2(\Psi_s(\mu),\Psi_t(\nu))^2
=s t W_2(\mu,\nu)^2 
+(s-t)
(s W_2(\delta,\mu)^2-t W_2(\delta,\nu)^2).
\]
\end{lemma}
\begin{proof}
Since the case that $s t$ equals $0$ is trivial, 
we consider the case that $st$ is positive.
It follows from the definition of the cone distance that
\begin{align*}
 d(\psi_s(x_1),\psi_t(x_2))^2 
&=s t d(x_1, x_2)^2
+(s-t) \left( s d(o, x_1)^2- t d(o, x_2)^2 \right)
\end{align*}
for any $x_1, x_2$ in $X$.
Because 
$(\psi_s \times \psi_t)_{\sharp} \pi$ 
is a transport plan in $\Pi(\Psi_s(\mu),\Psi_t(\nu))$
for any optimal transport $\pi$ in $\Pi(\mu,\nu)$, 
we obtain
\begin{align}\label{hutousiki} \notag
 \  W_2(\Psi_s(\mu),\Psi_t(\nu))^2  \notag
 &\leq \int_{X \times X} 
  d(x_1,x_2)^2 d (\psi_s \times \psi_t)_{\sharp} \pi(x_1,x_2)     \\ \notag
 &= \int_{X \times X} 
  d(\psi_s(x_1),\psi_t(x_2))^2 
  d \pi (x_1,x_2)\\
 &= s t W_2(\mu,\nu)^2+(s-t)(s W_2(\delta,\mu)^2-t W_2(\delta,\nu)^2).
\end{align}
The last equality follows from the facts
that $\pi$ is optimal and 
that the marginals of 
$\pi$ are $\mu$ and $\nu$.
We also obtain,
by substituting $(1/s,1/t)$ for $(s,t)$ and $(\Psi_s(\mu), \Psi_t(\nu))$ for $(\mu, \nu)$
respectively in \eqref{hutousiki}, that
\begin{align}\label{ineq}
 W_2(\mu,\nu)^2 
 \leq \frac{1}{s t} W_2(\Psi_s(\mu),\Psi_t(\nu))^2
+\left(\frac1s-\frac1t\right)(s W_2(\delta,\mu)^2
  -t W_2(\delta,\nu)^2).
\end{align}
Combining \eqref{hutousiki} and \eqref{ineq}, we deduce 
\begin{align*}
 W_2(\Psi_s(\mu),\Psi_t(\nu))^2 
 &\leq s t W_2(\mu,\nu)^2 
 +(s-t)
 (s W_2(\delta,\mu)^2-t W_2(\delta,\nu)^2) 
 \leq W_2(\Psi_s(\mu),\Psi_t(\nu))^2 . 
\end{align*}
Therefore the previous inequalities have to be equalities. 
\end{proof}
\begin{proof}[of the ``if'' part of Theorem~\ref{ippan}]
First of all, Lemma~\ref{kyori} guarantees 
the uniqueness of geodesics connecting $\delta$ and 
any $\mu$ in $\px$.
To see this, let $\{\mu(s)\}_{s \in [0,1]}$ be 
a geodesic from $\delta$ to $\mu$.
Then we have 
$W_2(\Psi_s(\mu),\mu(s))^2=0$
because $\mu=\Psi_1(\mu)$.
Thus all geodesics starting at $\delta$ 
are written as $\{\Psi_s(\mu)\}$ and
they can be extended up to the boundary of 
the ball $B(\delta,R)$ for any positive number $R$.
We may without loss of generality 
choose a geodesic ray $\gamma_{\mu} (s)=\Psi_s(\mu)$ 
from $\delta$ passing through $\mu$ with $W_2(\delta,\mu)=1$ 
as a representative of the equivalence classes 
in $\Sigma^{\prime}_{\delta}$.
We moreover conclude that 
the angle between $\gamma_{\mu}$ and $\gamma_{\nu}$ 
is given by
\begin{equation}\label{angle}
 \angle_{\delta} (\gamma_{\mu},\gamma_{\nu})=
 \arccos \left(1-\frac12W_2(\mu,\nu)^2 \right).
\end{equation}
The completeness of $(\px,W_2)$ 
guarantees the completeness of $(\Sigma^{\prime}_{\delta}, \angle_\delta)$ 
and it yields that 
the space of directions $\Sigma_{\delta}$ at $\delta$ 
in $\px$ coincides with $\Sigma^{\prime}_{\delta}$,
namely, $\Sigma_{\delta}$ is regarded as 
\[
  \left\{ \gamma_{\mu} \ | \ 
         \text{a geodesic ray from $\delta$ 
                passing through $\mu$ in $\px$ with 
               $W_2(\delta,\mu)=1$} 
 \right\}.
\]

We finally construct an isometric map $\Gamma$ 
from the tangent cone $(K_{\delta}, d_\delta)$ at $\delta$ to $(\px, W_2)$. 
Let $\Gamma$ be the map
given by $\Gamma(\gamma_{\mu},s)=\gamma_{\mu}(s)=\Psi_s(\mu)$,
which is well-defined and bijective. 
By Lemma~\ref{kyori} and \eqref{angle},
we get
\begin{align*}
  W_2(\Gamma(\gamma_{\mu},s),\Gamma(\gamma_{\nu},t))
 =W_2(\Psi_s(\mu),\Psi_t(\nu))
 =d_{\delta}((\gamma_{\mu},s),(\gamma_{\nu},t)),
\end{align*}
proving that $\Gamma$ is an isometry and 
the ``if'' part of Theorem~\ref{ippan}.
\end{proof}
%

When we prove the ``only if'' part of Theorem~\ref{ippan}, 
the following two lemmas play essential roles.
Although they are special cases of~\cite[Lemma~2.11]{st},
we include the proofs of them for the completeness of the argument.
\begin{lemma}\label{endpoints}
Let $(X, d)$ be a Polish space.
For arbitrary points $x$ and $\x$ in $X$, 
if there exists a unique measure $\mu$ in $\px$ so that 
\begin{align*}
 W_2(\mu,\delta_{x})
=W_2(\mu,\delta_{\x})
=\frac12 W_2(\delta_{x}, \delta_{\x})
=\frac12 d(x,\x),
\end{align*}
then $\mu$ is a Dirac measure. 
\end{lemma}
\begin{proof}
By the direct calculation, we get 
\begin{align*}
 W_2(\delta_{x}, \delta_{\x})^2
&=2 W_2(\mu,\delta_{x})^2 + 2 W_2(\mu,\delta_{\x})^2\\
&=2 \int_{X} 
             d(x,y)^2 d\mu(y)
 +2 \int_{X} 
             d(y,\x)^2  d\mu(y)\\
&\geq \int_{X} 
             d(x,\x)^2  d\mu(y)\\
&= d(x,\x)^2. 
\end{align*}
The inequality must be equality and we obtain 
\begin{align*}
 d(x,y)=d(y,\x)
=\frac12 d(x,\x),
\end{align*}
for $\mu$-almost every $y$ in $X$.
By the uniqueness of $\mu$, 
$\mu$ must be the Dirac measure. 
\end{proof}

\begin{lemma}\label{support}
Let $(X, d)$ be a Polish space and 
$\{\mu(s)\}_{s \in [0,1]}$ be a geodesic in $(\px, W_2)$.
If the midpoint $\mu(1/2)$ is a Dirac measure $\delta_{\x}$,
then for every $x_i$ in the support of $\mu(i)$ $(i=0,1)$,
we have 
\begin{align*}
W_2(\mu(i),\delta_{\x})
=d(x_i,\x)
=\frac12 d(x_0,x_1)
=\frac12 W_2(\mu(0),\mu(1)).
\end{align*}
In particular, if $X$ is non-branching at $\x$ 
in the sense of \eqref{non-branch}, 
then $\mu(0)$ and $\mu(1)$ are also Dirac measures.
\end{lemma}
\begin{proof}
Since the product measure $\mu=\mu(0)\times \mu(1)$ 
is a transport plan in $\Pi(\mu(0),\mu(1))$,
we get the following inequalities: 
\begin{align*}
 W_2(\mu(0),\mu(1))
&\leq \left( \int_{X \times X} 
             d(x_0,x_1)^2 d\mu(x_0,x_1)
      \right)^{\frac12} \\
&\leq \left( \int_{X \times X} 
             d(\x,x_0)^2  d\mu(x_0,x_1)
       \right)^{\frac12} 
 +\left( \int_{X \times X} 
             d(\x,x_1)^2  d\mu(x_0,x_1)
       \right)^{\frac12}\\
&= W_2(\mu(0),\mu(1/2))+W_2(\mu(1/2),\mu(1))\\
&=W_2(\mu(0),\mu(1)).
\end{align*}
The previous inequalities must be equalities and 
since the integrands are continuous,
we obtain that
\begin{equation*}\label{midpt}
d(x_0,\x) = d(\x,x_1) = \frac12 d(x_0,x_1)
\end{equation*}
for every $(x_0,x_1)$ in the support of $\mu= \mu(0) \times \mu(1)$.
It follows from this that
$d(x_0,x_1)$ is independent of $x_0$ and $x_1$,
and equal to $W_2(\mu(0),\mu(1))$.
This is the desired result.
\end{proof}

\begin{corollary}\label{vertex}
For a Polish space $(X, d)$, 
if its $L^2$-Wasserstein space has a cone structure 
then the element in $\px$ corresponding to the vertex must be a Dirac measure.
\end{corollary}
\begin{proof}
For an arbitrary point $x$ in the support of 
the vertex $\mu$, 
there exists a geodesic ray from $\mu$ 
passing through the Dirac measure $\delta_x$ in $\px$.
By applying Lemma~\ref{support}, we acquire 
\[
 W_2(\mu,\delta_x)=d(x,x)=0.
\]
Therefore the vertex is a Dirac measure. 
\end{proof}
\begin{proof}[of the ``only if'' part of Theorem~\ref{ippan}]
We assume that $(\px,W_2)$ is isometric to a cone
$(C(\Sigma),d_C)$ with the vertex $O$, 
which corresponds to a Dirac measure $\delta$ 
by Corollary~\ref{vertex}.
For any $x$ in $X$, 
there exists $(\xi, s_0)$ in $C(\Sigma)$ 
corresponding to the Dirac measure $\delta_x$.
Let $\Gamma_{\xi} : [0,\infty) \to C(\Sigma)$ be 
the geodesic ray given by $\Gamma_{\xi}(s)=(\xi,s)$.
Then for any non-negative $s$,
$\Gamma_{\xi}(s)$ also corresponds to a Dirac measure in $\px$.
This is due to the uniqueness of geodesics
in the cone $C(\Sigma)$ starting at the vertex,
as well as Lemma~\ref{endpoints} for $s$ in $(0, s_0)$
and Lemma~\ref{support} for $s$ in $(s_0, \infty)$.

We set 
\begin{align*}
Y
&=\{\xi \in \Sigma \ |\ (\xi,1) 
    \text{ corresponds to some Dirac measure}  \}\\
&=\{\xi \in \Sigma \ |\ (\xi,s) 
    \text{ corresponds to some Dirac measure for all } s>0  \}
\end{align*}
and define a map $\Gamma$ from $(C(Y),d_C)$ to $(X, d)$ by 
$\Gamma(\xi,s)=\gamma_\xi(s)$, 
where $\gamma_\xi(s)$ is the center of the Dirac measure in $\px$
corresponding to $\Gamma_{\xi}(s)$.
(We abbreviate the distance $d_C|_{C(Y)}$ 
as $d_C$.)
By the previous argument, 
$\Gamma$ is well-defined and bijective.
We additionally have 
\begin{align*}
 d_C((\xi,s),(\eta,t))
=W_2(\delta_{\gamma_\xi(s)},\delta_{\gamma_\eta(t)}) 
=d(\gamma_\xi(s),\gamma_\eta(t)) 
=d(\Gamma(\xi,s),\Gamma(\eta,t))
\end{align*}
for any $(\xi,s)$ and $(\eta,t)$ in $C(Y)$.
This shows that $\Gamma$ is an isometry between $C(Y)$ and $X$.
Now the proof of the ``only if'' part of Theorem~\ref{ippan} is complete.

The second statement about being non-branching at the vertices 
is a immediate consequence of Lemmas~\ref{support} and~\ref{vertex}. 
\end{proof}
\begin{remark}
When the space $(X,d_X)$ is a cone over $(Y,d_Y)$,
a relation between 
$(\py, W_2^Y)$ and the space of directions 
$(\Sigma_{\delta}, \angle)$ which is the base space of $(\px,W_2^X)$ is as follows.
We denote by $o$ and $\delta=\delta_o$ 
the vertex of $X$ and the Dirac measure centered at $o$, 
respectively.
%

We define two maps $\iota$ and $\Theta$ as follows:
\[
 \iota : Y \ni y \mapsto (y,1) \in X=C(Y), \quad 
 \Theta : \py\ni\mt \mapsto 
           \gamma_{\iota_{\sharp}\mt} \in \Sigma_{\delta}.
\]
The map $\Theta$ is well-defined, that is,
$W_2^X(\delta,\iota_{\sharp}\mt)=1$ for all $\mt$ in $\py$, and injective.
%

The map $\Theta: \py \to \Sigma_\delta$ as above is not an isometry in general.
Indeed, for 
\[
 Y=\s^1 = [-\pi,\pi]/\{\pi =- \pi\}, \quad 0<\theta<\pi/3,
\]
we define probability measures
$\mt,\nt$ on $Y$ by
\begin{align*}
 \mt =\frac12 \left( \delta_{y_1} + \delta_{y_2} \right),\quad  
 \nt =\frac12\left(\delta_{\y_1}+\delta_{\y_2}\right),
\end{align*}
where 
$(y_1,y_2)=(0,\pi-2\theta)$ and
$(\y_1,\y_2)=(\theta,\pi)$. 
Then we have 
\[
 W_2^Y(\mt,\nt)^2
=\frac52 \theta^2
\quad\text{ and }\quad 
 \cos\angle(\Theta(\mt),\Theta(\nt))
=\frac12(\cos \theta+\cos 2\theta),
\]
that is, 
 $\angle(\Theta(\mt),\Theta(\nt))$ 
is not equal to  $W_2^Y(\mt,\nt)$ and
$\Theta$ does not have a monotonicity of distance;
$W_2^Y(\mt, \nt)$ is smaller than $\angle(\Theta(\mt), \Theta(\nt))$ for $\theta$ close to $\pi/3$,
while the reverse inequality holds if $\theta$ is a sufficiently small.
\end{remark}
\begin{remark}
Since the supports of elements in $\Theta(\py)$ 
are contained 
in $\partial B(o,1) \subset X$,
the map $\Theta$ in the previous remark is not surjective. 
To see this,
let $N$ be the normal distribution on $\R$,
namely, its Radon--Nikodym derivative 
with respect to the Lebesgue measure $d t$ is given by
\[
 \frac{d N}{d t}(t)
=\frac1{\sqrt{2\pi}}\exp\left(-\frac{t^2}{2}\right)
\] 
and 
$\tau$ be the map from $\R$ to $X$ 
sending $t$ to $(y,|t|)$ for some $y$ in $Y$.  
Then the push-forward measure $\tau_{\sharp}N$ 
belongs to $\Sigma_\delta=\partial B(\delta_o,1)$, 
not to $\Theta(\py)$.
\end{remark}
\begin{remark}
The cone structure has such a function as an inner product of a Hilbert space and
we are taking advantage of its $L^2$-structure characteristic 
in the proof of Lemma~\ref{kyori}.
Thereby we do not expect to generalize our result 
to $L^p$-Wasserstein spaces.
\end{remark}
\section{Applications}
We first consider the next lemma 
which is the key of the proof of Corollary~\ref{rank}.
\begin{lemma}\label{ranks}
For any isometric embedding $\varphi$
of a cone $(C(Z), d_{C(Z)})$ over $(Z,d_Z)$ into 
another cone $(C(Y), d_{C(Y)})$ over a complete metric space $(Y,d_Y)$, 
there exists an isometric embedding $\psi$
of $C(Z)$ into $C(Y)$ 
which maps the vertex of $C(Z)$ to the vertex of $C(Y)$.
\end{lemma}
\begin{proof}
For any $(z, s)$ in $C(Z)$, 
we denote by $(\xi_z(s),r_z(s))$ its image $\varphi(z,s)$ in $C(Y)$.
Due to the triangle inequality, 
we acquire
\[
 |r_z(s)-s|\leq r_z(0)
\quad\text{ and hence }\quad
\lim_{s\to \infty }\frac{r_z(s)}{s}=1.
\]
Since
\begin{align*}
\cos (\min \{ d_Y(\xi_z(i),\xi_z(j)),\pi\})
&= \frac{r_z(i)^2 + r_z(j)^2 - (i-j)^2}{2r_z(i)r_z(j)}\\
&= \frac{ r_z(i)^2-i^2 + r_z(j)^2-j^2 +2ij}{2r_z(i)r_z(j)}
\to 1
\end{align*}
as $i$ and $j$ tend to infinity,
we notice that
$\{\xi_z(i)\}_{i\in \N}$ is a Cauchy sequence in $Y$.
Then the completeness of $Y$ ensures
the existence of the limit $\xi_{z}$ 
of $\{\xi_z(i)\}_{i\in \N}$.

Defining the map $\psi$ from $C(Z)$ to $C(Y)$ by $\psi(z,s)=(\xi_z, s)$,
we obtain the following equalities:
\begin{align*}
d_{C(Y)}(\psi(z,s),\psi(\zp,\tp))
&=s^2+\tp^2- 2s\tp 
 \lim_{i \to \infty} \cos (\min \{ d_Y(\xi_z(i),\xi_{\zp}(i)),\pi\})\\
&=s^2+\tp^2- 2s\tp\cos (\min \{ d_Z(z,\zp),\pi\})\\
&=d_{C(Z)}((z,s),(\zp,\tp)).
\end{align*}
This implies that the map $\psi$ is the desired isometric embedding 
which sends the vertex of $C(Z)$ to that of $C(Y)$.
\end{proof}
\begin{proof}[of Corollary~\ref{rank}]
It suffices to consider the ranks of $X$ and $\px$.
For any non-negative integer $k$ less than
or equal to the rank of $\px$, 
there exists an isometric embedding 
of $\R^k$ into $\px$ sending 
the origin $0$ to the vertex 
by the definition of the rank and Lemma~\ref{ranks}. 
In addition, Lemma~\ref{support} asserts that 
all of the elements lying in the image of $\R^k$ are Dirac measures. 
Thus the rank of $X$ is larger than or equal to $k$,
proving the equality of the ranks of $X$ and $\px$.
\end{proof}

We next prove Theorem~\ref{split}, 
which is a generalization of~\cite[Lemma~5(b)]{cmv}.
\begin{proof}[of Theorem~\ref{split}]
Changing notation, we assume that 
$(X, d_X)$ is isometric to the direct product of 
some metric space $(Y,d_Y)$ and $(\h, \lr)$.
For any $x$ in $X$, 
let $x_Y$ and $x_{\h}$ stand 
for the projections of $x$ to $Y$ and $\h$, respectively.
Then for any $x,\x$ in $X$, 
we have 
\[
 d_X(x,\x)^2
=d_Y(x_Y,\x_Y)^2 
+\|x_\h -\x_\h\|^2,
\]
where, $\|\cdot \|$ is the norm of $\h$.
Since $\h$ is the Hilbert space, 
an arbitrary $\mu$ in $\px$ has 
a unique mean $m(\mu)$ in $\h$ satisfying
\[
 \langle m(\mu),h \rangle
=\int_{X} \langle x_\h, h \rangle d \mu(x)
\]
for any $h$ in $\h$. 
We denote by $\ppz$ the subset of $\px$
of elements whose means are $h$.
We define a map $\varphi_h$ on $X$ 
by $\varphi_h(x)=(x_Y, x_{\h}+h)$ 
and the associated map $\Phi_{h}$ 
from $\px$ to $\px$ by $\Phi_{h}(\mu)=\vzs \mu$.
Then we have for any $\z$ in $\h$
\begin{align*}
 \int_{X} \langle x_\h,\z \rangle d\Phi_h(\mu)(x)
 =\int_{X} \langle x_\h+h,\z \rangle d\mu(x) 
 =\langle m(\mu)+h,\z \rangle,
\end{align*}
implying $m(\Phi_h(\mu))= m(\mu)+h$.
Since the maps $\Phi_h$ and $\Phi_{-h}$ are inverses of each other,
we conclude
\begin{align*}
 \px
 =\bigsqcup_{h\in\h} \ppz 
 =\bigsqcup_{h\in\h} \Phi_h(\ppo), 
\end{align*}
where $0$ is the zero vector in $\h$.
Thus we can define a map $\Phi$ from $\ppo \times \h$ to $\px$ 
by sending $(\mu,h)$ to $\Phi_h(\mu)=\mu_h$.
Now we confirm that the map $\Phi$ is an isometry.
For any $\mu$ and $\mmp$ in $\ppo$, 
$(\vz\times\vo)_{\sharp}\pi$ is a 
transport plan in $\Pi(\mu_h,\mmp_{\z})$ 
for any optimal transport plan $\pi$ 
in $\Pi(\mu,\mmp)$.
Then we get 
\begin{align*}
 W_2(\mu_h,\mmp_\z)^2 
 &\leq  \int_{X \times X} d_X(x,\x)^2 
         d (\vz\times\vo)_{\sharp}\pi  (x,\x)\\
 &= \int_{X \times X} 
   \Bigl[d_Y(x_Y,\x_Y)^2 
   +\|(x_\h+h) -(\x_\h+\z)\|^2 \Bigr]  d\pi (x,\x)\\
 &= \int_{X \times X}\Bigl[
     d_Y(x_Y,\x_Y)^2 
     +\|x_\h-\x_\h\|^2 
     +\|h-\z\|^2
     +2\langle x_\h-\x_\h,h-\z \rangle  \Bigr] 
       d\pi (x,\x)\\
 &= \int_{X \times X} 
     d_X(x,\x)^2 d\pi (x,\x) 
     +\|h-\z\|^2  \\
 &= W_2(\mu,\mmp)^2 +\|h-\z\|^2.
\end{align*}
The third equality follows from 
$m(\mu)=m(\mmp)=0$ and 
the last equality follows from the fact
that $\pi$ is optimal.
By a similar argument,
we also obtain 
\begin{align*}
  W_2(\mu,\mmp)^2 
\leq W_2(\mu_h,\mmp_\z)^2 -\|h-\z\|^2.
\end{align*}
Therefore we acquire 
\begin{align*}
  W_2(\mu_h,\mmp_\z)^2 
 \leq W_2(\mu,\mmp)^2 +\|h-\z\|^2 
 \leq W_2(\mu_h,\mmp_\z)^2.
\end{align*}
Hence the previous inequalities have to be equalities,
that is, the map $\Phi$ is an isometry.
\end{proof}
\begin{remark}\label{cg}
Carlen--Gango~\cite{cg} also investigate the structure of 
the absolutely continuous part of the $L^2$-Wasserstein spaces
over $\R^d$.
They do this in order to 
carry out the constrained version of the variational scheme of
Jordan--Kinderlehrer--Otto~\cite{jko}.

Our result seems to have in common with
what was established in~\cite[Section~3]{cg} (cf.~\cite[p.219, Line~12]{cmv}).
However,
our argument, based on the metric geometry, is simple
and works as well for the measures not necessarily absolutely continuous
with respect to the Lebesgue measure.
It will also help the readers interested in the results of~\cite{cg}.

See also the subsequent paper~\cite{cg2}
and Tudorascu's paper~\cite{T} where some open problems in~\cite{cg} are solved.
\end{remark}

\begin{remark}
If we choose a separable Hilbert space $\h$ as $X$ in the statement of Theorem~\ref{split},
then it turns out 
that $\mathcal{P}_{2,0}(\h)$ has a cone structure.
Moreover the base space $\Sigma_0$ of $\mathcal{P}_{2,0}(\h)$ is given by
\[
 \{\text{$\gamma_{\mu}$} \ | \ 
 \text{a geodesic ray from $\delta=\delta_0$ 
      in $\ppo$ passing through $\mu$ with $W_2(\delta,\mu)=1$}  \}
\]
by a similar argument 
in the proof of ``if'' part of
Theorem~\ref{ippan}.
We can estimate the diameter of $(\Sigma_0, \angle)$.
This estimate in the case of $\R^d$ appears in~\cite[(3.10)]{cg}, 
however it was proved in a different way.
For any elements 
$\gamma_\mu$ and $\gamma_\mmp$ in $\Sigma_0$, 
we acquire 
\begin{align*}
 W_2(\mu,\mmp)^2 
\leq \int_{\h \times \h} \|h-\z\|^2 d(\mu \times \mmp)(h,\z)
=1+1-2 \int_{\h \times\h} \langle h, \z \rangle d(\mu \times \mmp)(h,\z)
=2
\end{align*}
because the means and variances of $\mu$ and $\mmp$
are $0$ and $1$, respectively.
By \eqref{angle}, we obtain 
\begin{align*}
\cos \angle (\gamma_\mu, \gamma_\mmp) 
=1-\frac12W_2(\mu,\mmp)^2 \geq 0.
\end{align*}
It implies that the angle 
$\angle (\gamma_\mu, \gamma_\mmp )$ 
is smaller than or equal to $\pi/2$.
Since $\mu$ and $\mmp$ are arbitrary,
we acquire 
\[
 \mathrm{diam}  \Sigma_0
=\sup_{\mu,\nu \in \Sigma_0} 
 \angle (\gamma_\mu, \gamma_\mmp ) \leq \frac{\pi}{2}
\]
This fact and Theorem~\ref{split} also guarantee 
that the rank of $\pp$ is equals to $d$.
\end{remark}

\begin{proof}[of Corollary~\ref{tilps}]
Since part~(i) follows immediately from
Lemma~\ref{support}, we only prove part~(ii).
We assume that 
$\px$ has a cone structure and is isometric to the direct product of 
some metric space $(\mathcal{Q},d_\mathcal{Q})$ and $\h$.
We denote by $(q_0,0)$ 
the element in $\mathcal{Q}\times \h$ corresponding to the vertex of $\px$.
For any $h$ in $\h$,
the map given by $(q,\z) \mapsto (q,h+\z)$ is an isometry;
therefore $(q_0,h)$ must correspond to the vertex, 
i.e., a Dirac measure (Corollary~\ref{vertex}).

Next,
we fix an arbitrary $x$ in $X$
and find $(q,h)$ in $\mathcal{Q}\times \h$ corresponding to $\delta_x$ in $\px$.
Then there exists 
a unique geodesic ray $\{(q(s),0)\}_{s \in [0, \infty)}$ 
from the vertex $(q_0,0)$ passing through $(q,0)$ 
at $s=s_0=d_\mathcal{Q}(q_0,q)$.
By the argument used in the proof of the ``only if'' part of Theorem~\ref{ippan},
we can conclude that 
$(q(s), \frac{s}{s_0} h)$ must be corresponding to a Dirac measure 
for any positive number $s$.
Then consider the geodesic ray
$\{\ell(s)\}_{s\in [0, \infty)}$ in $\mathcal{Q}\times \h$ given by
\[
 \ell(s)=\left( q(s),
          2 \left(\frac{s}{s_0}-1\right)h
              +\left(2-\frac{s}{s_0}\right)\z \right)
\]
for some fixed $\z$ in $\h$.
We deduce that 
$\ell(s)$ is also corresponding to a Dirac measure in $\px$
for any positive $s$, 
because
$\ell(0)=(q_0,2(\z-h))$ and $\ell(2s_0)=(q(2s_0),2h)$ 
correspond to Dirac measures.
Thus $\ell(s_0)=(q,\z)$ is corresponding to a Dirac measure and
$X$ is isometric to 
\[
 \h \times
\{q \in \mathcal{Q} \ |\  
 \text{$(q,h)$ corresponds to some Dirac measure 
        for some (hence all) $h \in \h$} \}.
\]
This completes the proof of part~(ii) of Corollary~\ref{tilps}.
\end{proof}

\begin{remark}
We mention the result of Mitsuishi~\cite{Mi}, 
which is the splitting theorem for 
Alexandrov spaces of non-negative curvature 
without the properness assumption
(see~\cite{Mi} for the statement and definitions).
With the help of Corollary~\ref{tilps}.(i),
his result ensures that  
the rank of the $L^2$-Wasserstein space over any Alexandrov space of non-negative curvature
is the same as that of the underlying space,
since being non-branching everywhere is one of the fundamental properties of 
Alexandrov spaces with lower curvature bound.

\end{remark}

\def\cprime{$'$}
\providecommand{\bysame}{\leavevmode\hbox to3em{\hrulefill}\thinspace}
\providecommand{\MR}{\relax\ifhmode\unskip\space\fi MR }
\providecommand{\MRhref}[2]{%
  \href{http://www.ams.org/mathscinet-getitem?mr=#1}{#2}
}
\providecommand{\href}[2]{#2}

\affiliationone{
   A. Takatsu\\
  Mathematical Institute\\
  Tohoku University\\ 
  Sendai 980-8578\\
  Japan
   \email{sa6m21@math.tohoku.ac.jp}}
\affiliationtwo{
   T. Yokota\\
   Graduate School of Pure and Applied Sciences\\
   University of Tsukuba\\
   Tsukuba 305-8571\\
  Japan
   \email{takumiy@math.tsukuba.ac.jp}}
\end{document}